\input amstex
\input epsf
\input Amstex-document.sty

\pageno 775

\def\shorttitle{Subfactors and Planar Algebras}
\def\shortauthor{D. Bisch}

\topmatter %
\title\nofrills{\boldHuge Subfactors and Planar Algebras}
\endtitle

\author \Large D. Bisch* \endauthor

\thanks *Department of Mathematics, Vanderbilt University,
Nashville, TN 37240 and UCSB, Department of Mathematics, Santa Barbara, CA 93106, USA. E-mail:
bisch\@math.vanderbilt.edu \endthanks

\abstract\nofrills \centerline{\boldnormal Abstract}

\vskip 4.5mm

An inclusion of II$_1$ factors $N \subset M$ with finite Jones index
gives rise to a powerful set of invariants that can be approached
successfully in a number of different ways. We describe Jones' pictorial
description of the standard invariant of a subfactor as a so-called planar
algebra and show how this point of view leads to new structure
results for subfactors.

\vskip 4.5mm


\endabstract
\endtopmatter

\document

\baselineskip 4.5mm \parindent 8mm

\specialhead \noindent \boldLARGE 1. Introduction \endspecialhead

{\it Abelian} von Neumann algebras are simply algebras of bounded,
measurable functions on a measure space.
A general (non-abelian) von Neumann algebra can be viewed as an algebra
of ``functions'' (operators) on a {\it non-commutative measure space}.
The building blocks of what one might call
{\it non-commutative probability spaces} are the so-called
II$_1$ {\it factors} $M$,
that is those von Neumann algebras with trivial center that are infinite
dimensional and possess a distinguished tracial state (the analogue
of a {\it non-commutative integral}). The ``smallest'' II$_1$ factor
is the {\it hyperfinite} II$_1$ factor which is obtained as
the closure in the weak operator topology of the canonical
anti-commutation relations (CAR) algebra
of quantum field theory. A II$_1$ factor comes always with a natural
left representation on $L^2(M)$, the non-commutative L$^2$-space associated
to $M$. See for instance [13].

Vaughan Jones initiated in the early 80's the {\it theory of subfactors}
as a ``Galois theory'' for inclusions of II$_1$ factors. A {\it subfactor}
is an inclusion of II$_1$ factors $N \subset M$
such that the dimension of $M$ as left $N$-Hilbert module is finite.
This dimension is called the {\it Jones index} $[M:N]$ ([19]) and
one would expect by classical results of Murray and von Neumann that it takes
on any real number $\ge 1$. One of the early results in
the theory of subfactors was Jones' spectacular {\it rigidity theorem} which
says that this index is in fact {\it quantized} [19]:
if $[M:N] \le 4$, then it has to be of the form $4 \cos^2 \frac{\pi}{n}$,
for some $n \ge 3$.
Since Jones' early work the theory of subfactors has developed into one of
the most exciting and rapidly
evolving areas of operator algebras with numerous applications to different
areas of mathematics (e.g. knot theory with the discovery of the
{\it Jones polynomial} [20]), quantum physics and statistical mechanics.
Subfactors with finite
Jones index have an amazingly rich mathematical structure and an interplay
of analytical, algebraic-combinatorial and topological techniques
is intrinsic to the theory.

\specialhead \noindent \boldLARGE 2. Subfactors \endspecialhead

A subfactor can be viewed as a group-like object that encodes what
one might call {\it generalized symmetries} of the data that went
into its construction. To decode this information
one needs to compute the {\it higher relative commutants}, a system
of inclusions of certain finite dimensional C$^*$-algebras naturally
associated to the subfactor. This system
is an invariant of the subfactor, the so-called {\it standard invariant},
which contains in many natural situations precisely the same
information as the subfactor itself ([30], [32], [33]). Here is one way to
construct the standard
invariant: If $N \subset M$ denotes an inclusion of II$_1$ factors
with finite Jones index, and $e_1$ is the orthogonal projection
$L^2(M) \to L^2(N)$, then we define $M_1$ to be the von Neumann algebra
generated by $M$ and $e_1$ on $L^2(M)$. $M_1$ is again
a II$_1$ factor and $M \subset M_1$ has finite Jones index as well so
that the previous construction can be repeated and {\it iterated} [19].
One obtains a tower of II$_1$ factors
$N \subset M \subset M_1 \subset M_2 \subset \dots$
associated to $N \subset M$, together with a remarkable sequence
of projections $(e_i)_{i \ge 1}$, the so-called {\it Jones projections},
which satisfy the Temperley-Lieb relations and give rise to
Jones' braid group representation [19], [20]. The
(trace preserving) isomorphism class of the system of inclusions
of (automatically finite dimensional)
centralizer algebras or {\it higher relative commutants}

\vskip -0.1in
$$\matrix {\Bbb C} = N' \cap N & \subset & N' \cap M & \subset &
N ' \cap M_1 & \subset & N' \cap M_2 & \subset & \cdots \\
&& \cup && \cup && \cup & \\
 && {\Bbb C} = M' \cap M & \subset & M' \cap M_1 & \subset &
M' \cap M_2 & \subset & \cdots \endmatrix $$
is then the standard invariant ${\Cal G}_{N,M}$ of the subfactor $N
\subset M$. Each row of inclusions is given by a sequence of Bratteli diagrams, which can in fact be reconstructed
from a single, possibly infinite, bipartite graph. Hence one obtains two graphs (one for each row), the so-called
{\it principal graphs} of $N \subset M$, which capture the inclusion structure of the above double-tower of higher
relative commutants. It turns out that if $M$ is hyperfinite and $N \subset M$ has {\it finite depth} (i.e. the
principal graphs are finite graphs) [30], [32] or more generally if $N \subset M$ is {\it amenable} [33], then the
standard invariant determines the subfactor. In this case the subfactor can be reconstructed from the finite
dimensional data given by ${\Cal G}_{N,M}$. In particular, subfactors of the hyperfinite II$_1$ factor $R$ with
index $\le 4$ are completely classified by their standard invariant and an explicit list can be given (see for
instance [14], [16] or [33]). If the Jones index becomes $\ge 6$ such an explicit list is out of reach as the work
in [6], [11] and [12] shows: there are uncountably many non-isomorphic, irreducible infinite depth subfactors of
$R$ with Jones index $6$ and the same standard invariant! Partial lists of irreducbile subfactors with index
between $4$ and $6$ have been obtained by different methods (see for instance [1], [5], [6], [17], [35], [36],
[37], [38]), but much work remains to be done.

There are several distinct ways to analyze the standard invariant of a
subfactor (see [2], [4], [14], [22], [30], [33]). For instance,
in the bimodule approach ([13], [30], see also [4], [14], [18])
${\Cal G}_{N,M}$
is described as a {\it graded tensor category} of natural bimodules
associated to the subfactor. ${\Cal G}_{N,M}$ can thus be viewed as an
{\it abstract system of (quantum) symmetries} of the
mathematical or physical situation from which the subfactor was constructed.
It is in fact a mathematical object which generalizes for instance
discrete groups and representation categories of quantum groups
([37], [38]). A variety of powerful and novel techniques have been
developed over the last years that make it possible to compute and understand
the standard invariant of a subfactor. A key result
is Popa's {\it abstract characterization} of the standard invariant [34].
Popa gives a set of axioms that an abstract system of inclusions of finite
dimensional C$^*$-algebras needs to satisfy in order to arise as
the standard invariant of some (not necessarily hyperfinite) subfactor.
This result makes it possible to analyze the structure of subfactors,
which are infinite dimensional, highly non-commutative objects,
by investigating the {\it finite dimensional} structures encoded in
their standard invariants.

\specialhead \noindent \boldLARGE 3. Planar algebras
\endspecialhead

Jones found in [22] a powerful formalism to handle complex computations
with ${\Cal G}_{N,M}$. He showed that the standard invariant of
a subfactor has an intrinsic {\it planar} structure (this will be made
precise below) and that certain {\it topological} arguments can be used
to manipulate the operators living in the higher relative commutants of
the subfactor. The standard invariant is a so-called
{\it planar algebra}. To explain this notion let us first define the
{\it planar ``operad''} following [22]. Elements of the planar operad
are certain classes
of {\it planar $k$-tangles} which determine multilinear operations
on the vector spaces underlying the higher relative commutants associated
to a finite index subfactor.

A planar $k$-tangle consists of the unit disk $D$ in the complex plane
together with several interior disks $D_1$, $D_2, \dots ,D_n$.
The boundary of $D$ is marked with $2k$ points and each $D_j$ has
$2k_j$ marked points on its boundary. These marked points are connected
by strings in $D$, which meet the boundary of each disk
transversally. We
also allow (finitely many) strings which are closed curves in the interior
of $D$. The main point is that all strings are required to be disjoint
(hence {\it planarity}) and to lie in the complement of the interiors
of the $D_j$'s. Additional data of a
planar $k$-tangle is a checkerboard shading of the connected
components of $ \overset \circ \to D \backslash \bigcup_{j=1}^{n} D_j$,
and a choice of a white region at every $D_j$ (which corresponds to
a choice of the {\it first} marked point on the boundary of each
$D_j$). The {\it planar operad} $\Cal P$ is defined to consist of all
orientation-preserving
diffeomorphism classes of planar $k$-tangles (for all $k \ge 0$), where
the diffeomorphisms leave the boundary of $D$ fixed but are allowed to
move the interior disks. $\Cal P$ becomes a {\it colored operad} [22] (see
[28]). An example of a $4$-tangle is depicted in the next figure:

\medskip
\centerline{ \epsfysize=120pt \epsffile{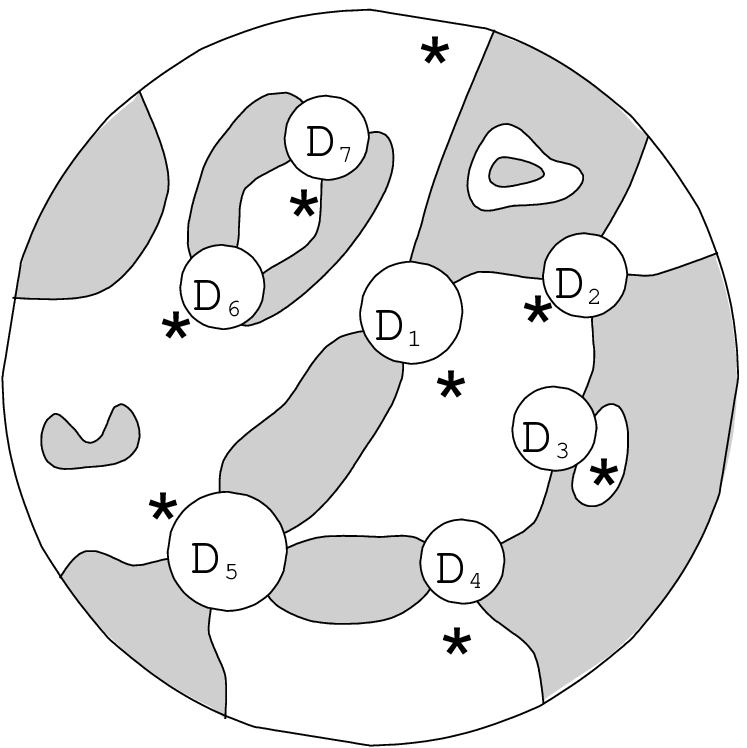}}
\medskip

Note that there are two classes of planar $0$-tangles according to the shading
of the tangle near the boundary of $D$.

Two planar tangles $\Cal T$ and $\Cal S$ can be composed in a natural
way if the number of boundary points of $\Cal S$ matches the
number of boundary points of one of the interior disks $D_j$ of $\Cal T$:
To obtain the composed tangle ${\Cal T} \circ_j {\Cal S}$
shrink $\Cal S$ and paste it inside $D_j$ so that the shadings and
marked white regions match up. Join the strings at the boundary of
$D_j$, smooth them and erase the boundary of $D_j$. It is clear
that this operation is well-defined (the checkerboard shading and choice
of a white region at each disk avoid rotational ambiguity) and that it
depends only on the
isotopy class of each tangle. Note that there may be several different ways
of composing two given tangles, each composition yielding potentially
distinct planar tangles. An example of such a composition is given in
the next figure (insert $\Cal S$ in the disk $D_2$ of $\Cal T$):

\medskip
\hskip -10pt{ \epsfysize=60pt \epsffile{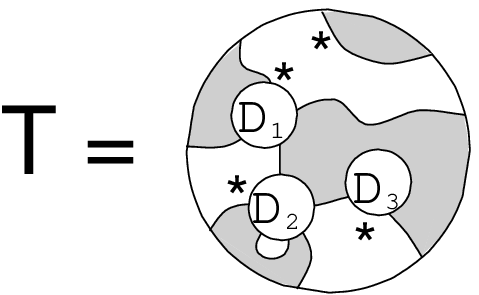} } \hskip 5pt{ \epsfysize=60pt \epsffile{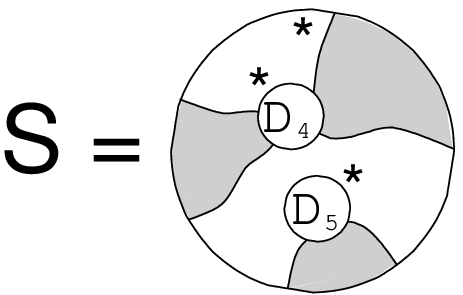} } \hskip
0.05in{ \epsfysize=70pt \epsffile{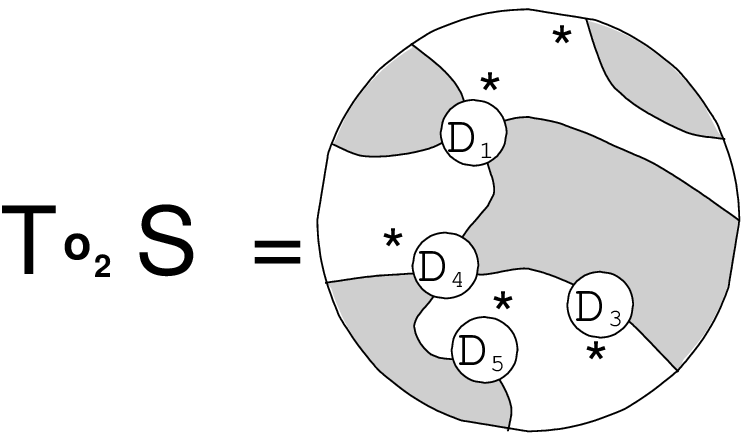} }

\medskip
An {\it abstract planar algebra} is then defined to be an algebra over this
planar operad ([28]). More concretely, an
abstract planar algebra $\Cal P$ is the disjoint union of vector
spaces ${\Cal P} = P_0^{\text \rm white} \coprod P_0^{\text \rm black}
\coprod_{n >0} P_n$
plus a morphism from the planar operad to the (colored) operad of
multilinear maps between these vector spaces. In other words
a planar algebra structure on $\Cal P$ is a procedure that assigns
to each planar $k$-tangle $\Cal T$ (with interior disks $D_j$ having
$2k_j$ boundary points, $1 \le  j \le n$) a multilinear map
$Z({\Cal T}): P_{k_1} \times \dots \times P_{k_n} \to P_k$ in such
a way that composition of tangles is compatible with
the usual composition of maps ({\it naturality} of composition).
Note that the $P_k$'s are
automatically associative algebras since the tangle in the next figure
(drawn in the case $k=5$)
defines an associative multiplication $P_k \times P_k \to P_k$
(associativity follows from naturality of the composition).

\medskip
\centerline { \epsfysize=80pt \epsffile{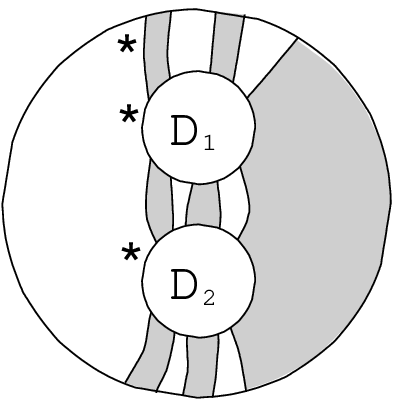} }


\medskip

Observe that this is a purely algebraic structure - the definition can be made for (possibly infinite dimensional)
vector spaces over an arbitrary field. The key point is of course that this structure appears naturally in the
theory of subfactors. In order to connect with subfactors several additional conditions will be required in the
definition of a planar algebra. A {\it planar algebra} (or {\it subfactor planar algebra} to emphasize the
operator algebra context) will be an abstract planar algebra such that $\dim P_k < \infty$ for all $k$, $\dim
P_0^{\text \rm white} = \dim P_0^{\text \rm black} = 1$ and such that the {\it partition function} $Z$ associated
to the planar algebra is positive and non-degenerate. The partition function is roughly obtained as follows: If
$\Cal T$ is a $0$-tangle, then $Z(\Cal T)$ is a {\it scalar} since it is an element in the 1-dimensional space
$P_0^{\text \rm white}$ resp. $P_0^{\text \rm black}$. Note that every planar algebra comes with two parameters
$\delta_1 = Z$(\epsfysize=10pt  \epsffile{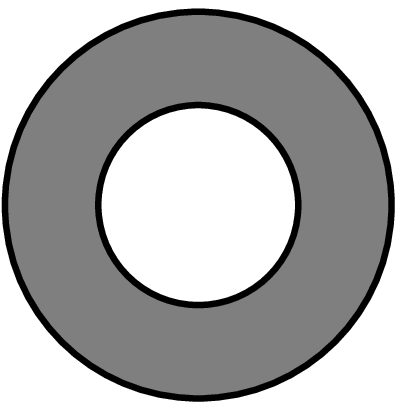}) and $\delta_2 = \break Z$(\epsfysize=10pt
\epsffile{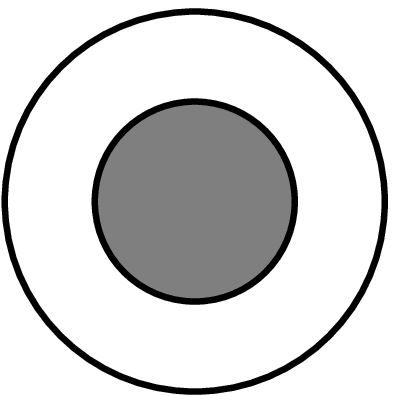}), which we require to be $\ne 0$ (the inner circles are strings, not boundaries of disks!).
In the case of a subfactor planar algebra we have $\delta \overset {\text \rm def} \to
 = \delta_1 = \delta_2$ (which is equivalent to extremality of
the subfactor [31]). In fact $\delta = [M:N]^{1/2}$ in this case. There
is an intrinsic way to define an involution on the planar algebra
arising from a subfactor which makes the partition function into
a sesquilinear form on the standard invariant. {\it Positivity of the
partition function} $Z$ means then positivity of this form. Note that
$Z$ gives in particular the natural trace on the standard invariant
of the subfactor. The main result of [22] is then the following theorem.

{\bf Theorem 3.1.} \it  The standard invariant ${\Cal G}_{N,M}$ of
an extremal subfactor $N \subset M$ is a subfactor planar algebra
$\Cal P = (P_n)_{n \ge 0}$ with $P_n = N' \cap M_{n-1}$.  \rm

This theorem says in particular that planar tangles {\it always} induce
multilinear maps (``planar operations'') on the standard invariant of
a subfactor.
As a consequence one obtains a diagrammatic formalism that can be
employed to manipulate the operators in $N' \cap M_{n-1}$ and
intricate calculations with these operators can be carried out using
simple topological arguments. This point of view has been turned in [9],
[10] into a powerful tool to prove general structure
theorems for subfactors, and to analyze the rather complex combinatorial
structure of the standard invariant of a subfactor. It has led to a
generators and relations approach to subfactors. See also [23], [24]
for more on this.

The two most fundamental examples of subfactor planar algebras are
the {\it Temperley-Lieb systems} of [19]
(see also [22]) and the {\it Fuss-Catalan systems} of [7] (see section 4).

Observe that by construction planar algebras are closely related to invariants
for graphs, knots and links and to the pictorial formalism commonly
used in the theory of integrable lattice models in statistical mechanics.

\specialhead \noindent \boldLARGE 4.  Fuss-Catalan algebras
\endspecialhead

Jones and I discovered in [7] a new hierarchy of
finite dimensional algebras, which arise as the higher relative
commutants of subfactors when intermediate subfactors are present.
These algebras have a number of interesting combinatorial properties
and they have recently been used to construct new integrable lattice models
and new solutions of the Yang-Baxter equation ([15], [29]).

We show in [7] that a chain of $k-1$ intermediate subfactors
$N \subset P_1 \subset P_2 \subset \dots P_{k-1} \subset M$ leads
to a tower of algebras $\big( FC_n(a_1,\dots , a_k) \big)_{n \ge 0}$,
which depend on $k$ complex parameters $a_1, \dots , a_k$.
The dimensions of these algebras are given by the generalized
Catalan numbers or {\it Fuss-Catalan numbers}
$\frac{1}{kn+1}\binom{(k+1)n}{n}$ and we therefore call these algebras
the {\it Fuss-Catalan algebras}.
If no intermediate subfactor is present, i.e. $P_i =N$ or $P_i =M$ for
all $i$, then one finds the well-known {\it Temperley-Lieb algebras}
(case $k=1$) [19]. The additional
symmetry coming from the intermediate subfactor is captured
{\it completely} by these new algebras and it is proved in
[7] (see also [8]) that they
constitute the minimal symmetry present whenever an intermediate
subfactor occurs. See also [26].

Let us explain in more detail what happens in the case of just one
intermediate subfactor. We consider $N \subset P \subset M$, an
inclusion of II$_1$ factors
with finite Jones index, and construct the associated tower of
of II$_1$ factors as in section 2. One obtains an inclusion of
II$_1$ factors
$N \subset P \subset M \overset p_1 \to \subset P_1 \overset e_1
\to \subset M_1  \overset p_2 \to \subset P_2 \overset e_2 \to \subset M_2
\subset \dots $, where the $p_i$'s are the orthogonal projections
from $L^2(M_{i-1})$ onto $L^2(P_{i-1})$ ($P_0=P$, $M_0 = M$) and
the intermediate subfactors $P_i$ are the von Neumann algebras
generated by $M_{i-1}$ and $p_i$. The algebra
$\text{\rm IA}_n(\alpha, \beta) \overset {\text{\rm def}} \to =
\text{ \rm Alg}(1,e_1, \dots,
e_{n-1}, p_1, \dots , p_{n-1})$, generated by the $e_i$'s and the
$p_i$'s, is a subalgebra of $N' \cap M_{n-1}$. It can be shown to depend
only on the two indices $\alpha=[P:N]$ and $\beta = [M:P]$, and {\it not}
on the particular position of $P$ in $N \subset M$. The projections $e_i$ and
$p_j$ satisfy again some rather nice commutation relations
(see [7] for details). In order to describe the structure
of these algebras let us for the moment consider the complex vector space
$FC_n(a,b)$, spanned by labelled, planar diagrams of the form

\medskip
\centerline{\eightpoint 2n marked points}

\vskip -0.35in \centerline{\epsfysize=60pt \epsffile{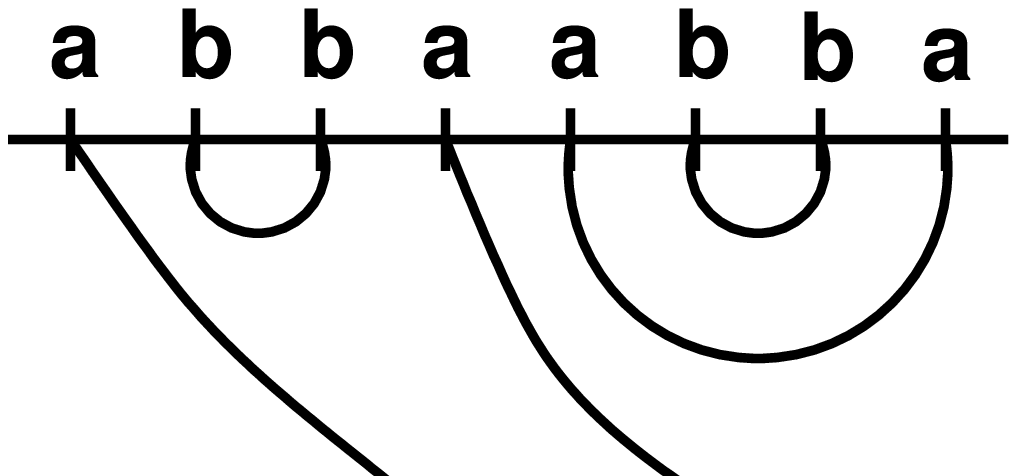}}
\bigskip
\bigskip
\bigskip

\noindent where $a$, $b \in \Bbb C \backslash \{ 0 \}$ are fixed. There is a natural multiplication of these
diagrams, which makes $FC_n(a,b)$ into an associative algebra (see [25]). To obtain $D_1 \cdot D_2$ put the basis
diagram $D_1$ on top of $D_2$ so that the labelling matches, remove the middle bar and all closed loops. Multiply
the resulting diagram with factors of $a$ resp. $b$ according to the number of removed $a$-loops resp. $b$-loops.
An example is depicted in the next figure.

\medskip
\centerline{{\epsfysize=90pt \epsffile{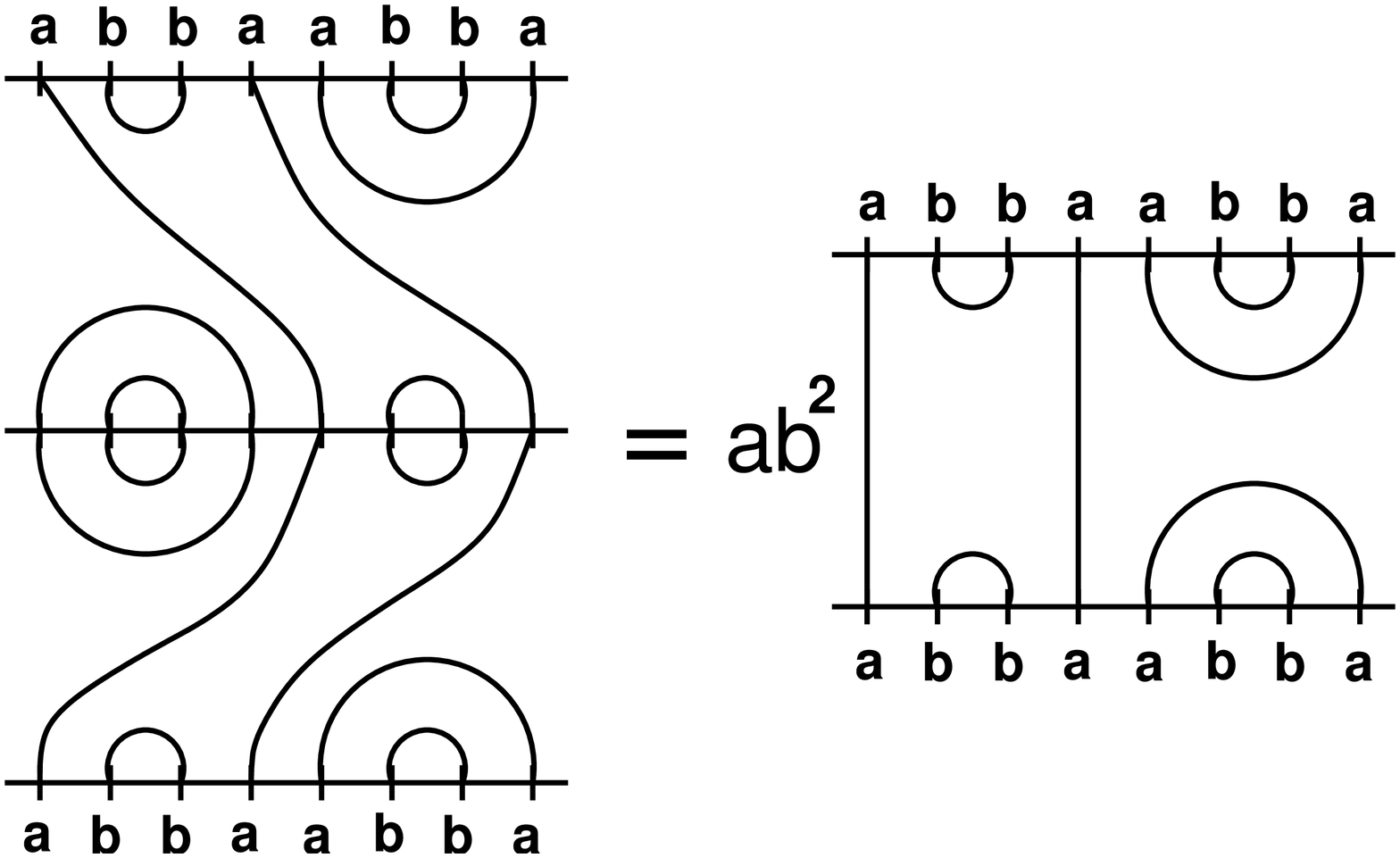}}}
\medskip

Counting diagrams shows that $\dim FC_n(a,b)= \frac{1}{2n+1}\binom {3n}n$, the n-th Fuss-Catalan number [7].
Clearly $FC_n(a,b)$ embeds as a subalgebra of $FC_{n+1}(a,b)$ by adding two vertical through strings to the right
of each basis diagram of $FC_n(a,b)$. A diagrammatic technique, called the {\it middle pattern analysis} in [7],
can be used to compute the structure of these algebras completely in the semi-simple case. One obtains that the
structure of the tower $FC_1(a,b) \subset FC_2(a,b) \subset \dots $ of Fuss-Catalan algebras is given by the {\it
Fibonacci graph} [7].

\medskip
\centerline{ {\epsfysize=130pt \epsffile{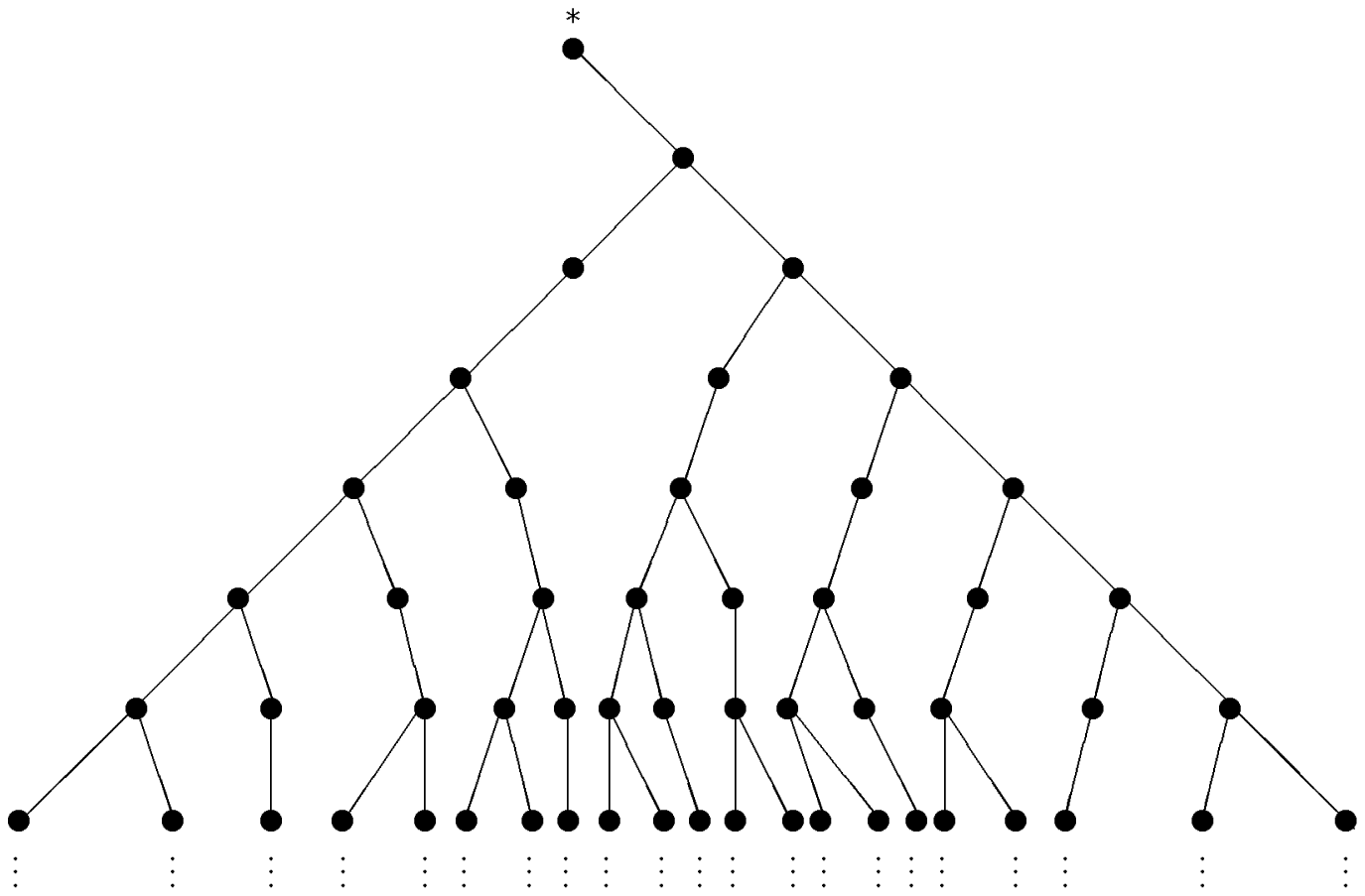}}}
\medskip

The algebras $\text{\rm IA}_n(\alpha, \beta)$ that we are interested in
can then be shown to be isomorphic to $FC_n(a,b)$, where $\alpha = {a^2}$,
$\beta = {b^2}$, if the indices $\alpha$ and $\beta$ are
generic, i.e. $> 4$. In the non-generic case $\text{\rm IA}_n(\alpha, \beta)$
is a certain quotient of $FC_n(a,b)$ (see [7] for the details).

There is a natural {\it 2-parameter Markov trace} on the Fuss-Catalan
algebras and the trace weights are calculated explicitly in [7].
In the special case of the Temperley-Lieb algebras this Markov trace
is the one discovered by Jones in [19]. The Fuss-Catalan tower together
with this Markov trace satisfies Popa's axioms in [34] and hence, one can
conclude from [34] that for every pair $(\alpha,\beta)$ of possible Jones
indices, there is a subfactor whose standard invariant is given precisely
by the corresponding Fuss-Catalan system
$\big(FC_n(\sqrt \alpha,\sqrt \beta ) \big)_{n \ge 0}$.
One obtains in this way uncountably many new subfactors.
A complete set of generators and relations for the Fuss-Catalan algebras
is also determined in [7].

It should be evident that the Fuss-Catalan algebras can be viewed
as planar algebras generated by a single element in $P_2 = N' \cap M_1$,
namely by the Jones projection $p_1$ onto the intermediate subfactor.
This projection can be characterized abstractly [3] and
it satisfies a remarkable exchange relation ([9], [27]), which
plays an important role in the work described in the next section.

\specialhead \noindent \boldLARGE 5. Singly generated planar algebras
\endspecialhead

Any subset $S$ of a planar algebra $\Cal P $ generates a planar subalgebra
as the smallest graded vector space containing $S$ and closed under planar
operations. From this point of view the simplest subfactors will
be those whose planar algebra is generated by the fewest elements
satisfying the simplest relations, while the index may be arbitrarily
large. If $S$ is empty we obtain the Temperley-Lieb algebra.
The next most complicated planar algebras after Temperley-Lieb should
be those generated by a single element $R$ which is in the
$k$-graded subspace $P_k$ for some $k>0$. We call such an element
a {\it $k$-box}. In [22] the planar algebra generated
by a single $1$-box was completely analyzed so the next case is
that of a planar algebra generated by a single $2$-box. This
means that the dimension of $P_2$ is at least $3$ so the first case
to try to understand is when $\dim P_3 = 3$. This dimension condition
by itself imposes many relations on $\Cal P$ but probably not enough to
make a complete enumeration a realistic goal. However, if one imposes
$\dim P_3 \leq 15$, then apart from a degenerate
case, this forces enough relations to reduce the number of variables
governing the planar algebra structure to be finite in number ([22],
see also [9]). It seems therefore reasonable to try to find
all subfactor planar algebras $\Cal P$ generated by a single element
in $P_2$ subject to the two restrictions
$\dim P_2 =3$ and $ \dim P_3 = d$ with $d \le 15$.
\medskip

In [9] we solved this problem when $d \leq 12$. In fact, using planar
algebra techniques we prove a much more general structure theorem for
subfactors.

{\bf Theorem 5.1.} \it
Let $N \subset M$ be an inclusion of II$_1$
factors with $3 < [M:N] < \infty$. Suppose that $\dim N' \cap M_1 = 3$
and that $N' \cap M_2$ is abelian modulo the basic construction
ideal $(N' \cap M_1)e_2(N' \cap M_1)$. Then there is an intermediate
subfactor $P$ of $N \subset M$, $P \ne N$, $M$. In particular
$Jx^*J=x$ for all $x \in N' \cap M_1$. \rm

The proof uses in a crucial way the abstract characterization
of the intermediate subfactor projection in [3] and
planar algebra techniques developed in [22] and [9].
It implies the following classification result.

{\bf Theorem 5.2.} \it If $\Cal P$ is a subfactor planar algebra
generated by a 3-dimensional $P_2$, subject to the condition
$\dim P_3 \le 12$, then it must be one of the following:
\roster
\item"a)" If $\dim P_3 = 9$, then it is the planar algebra associated
to the index $3$ subfactor $M^{\Bbb Z_3} \subset M$.
\item"b)" If $\dim P_3 =10$, then it is the $D_\infty$ planar algebra (a
special FC planar algebra).
\item"c)" If $\dim P_3 =11$ or $12$, then it is one of the FC planar algebras.
\endroster
\rm

The dimension conditions imply that a subfactor whose standard
invariant is a planar algebra of the form b) or c) satisfies the
hypothesis of Theorem 5.1 and hence must have an intermediate subfactor.
Since the Fuss-Catalan planar algebra is the minimal symmetry associated
to an intermediate subfactor it then follows easily that the planar
algebra has to be one of these.

It is quite natural to expect that increasing the dimension of $P_3$ should
result in a larger number of examples of planar algebras since there are
more a priori undetermined structure constants in the action of
planar tangles on $\Cal P$. Thus the result in [10] that there is a
{\it single} subfactor planar algebra satisfying the above restrictions
with $d=13$ is a complete surprise.
The planar algebra which arises is that of a subfactor obtained
as follows. Take an outer action of the dihedral group $D_5$ on a
type II$_1$ factor $R$ and let $M$ be the crossed product $R \rtimes D_5$ and
$N$ be the subfactor $R \rtimes \Bbb Z_2$. This particular subfactor has
played a significant role in the development of subfactors and relations with
knot theory and statistical mechanics. In [21] it was noted that there is a
solvable statisitical mechanical model associated with it and that it
corresponds to an evaluation of the Kauffman polynomial invariant of a link.
We prove in [10] the following

{\bf Theorem 5.3.} \it
Let ${\Cal P}= (P_k)_{k \ge 0}$ be a subfactor planar algebra generated
by a non-trivial element in $P_2$
(i.e. an element not contained in the Temperley-Lieb subalgebra of $P_2$)
subject to the conditions $\dim P_2 = 3$ and
$\dim P_3 = 13$. Then $\Cal P$ is the standard invariant of the
crossed product subfactor $R \rtimes {\Bbb Z}_2 \subset R \rtimes D_5$.
Thus there is precisely one subfactor planar algebra $\Cal P$ subject to
the above conditions. \rm

Note that this subfactor can be viewed as a Birman-Murakami-Wenzl
subfactor (associated to the quantum group of $Sp(4,\Bbb R)$ at a 5-th root
of unity, see [36]). We note here that the standard invariants
$\Cal P = (P_k)_{k \ge 0}$ of {\it all} BMW subfactors are generated
by a single non-trivial operator in $P_2$ and that they satisfy the
condition $\dim P_3 \le 15$.

The proof of this theorem uses in a crucial way theorem 5.1 and
the tight restrictions imposed by compatibility of
the rotation of period $3$ on $P_3$ and the algebra structure.

The next phase of this enumeration project will be to tackle the case
$d=14$. Here
we know that the quantum $Sp(4,\Bbb R)$ specialization of the BMW algebra
will give examples with a free parameter. We do expect however, that the
general ideas of [9] and [10] will
enable us to enumerate all such subfactor planar algebras.

\widestnumber\key{AAA}
\specialhead \noindent \boldLARGE References \endspecialhead

\ref
\key 1
\by M. Asaeda \& U. Haagerup
\paper Exotic subfactors of finite depth with Jones indices
$(5+\sqrt{13})/2$ and $(5+\sqrt{17})/2$
\jour Comm. Math. Phys.
\vol 202
\yr 1999
\pages 1--63
\endref

\ref
\key 2
\by T. Banica
\paper \rm Representations of compact quantum groups and
subfactors
\jour \it J. Reine Angew. Math.
\vol \rm 509
\yr 1999
\pages 167--198
\endref

\ref
\key 3
\by D. Bisch
\paper \rm A note on intermediate subfactors
\jour \it Pacific Journal of Math.
\pages 201-216
\vol \rm 163
\yr 1994
\endref

\ref \key 4 \by D. Bisch \paper \rm Bimodules, higher relative commutants and the fusion algebra associated to a
subfactor \inbook \it The Fields Institute for Research in Math. Sciences Commun. Series \vol \rm 13 \yr 1997,
13-63 \publ AMS, Providence, Rhode Island
\endref

\ref
\key 5
\by D. Bisch
\paper \rm An example of an irreducible subfactor of the hyperfinite
II$_1$ factor with rational, noninteger index
\jour \it J. Reine Angew. Math.
\vol \rm 455
\yr 1994
\pages 21-34
\endref

\ref
\key 6
\by D. Bisch \& U. Haagerup
\paper \rm Composition of subfactors: new examples of infinite depth subfactors
\jour \it Ann. scient. {\'E}c. Norm. Sup.
\vol \rm 29
\yr 1996
\pages 329-383
\endref

\ref
\key 7
\by D. Bisch \& V.F.R. Jones
\paper \rm Algebras associated to intermediate subfactors
\jour \it Invent. Math.
\vol \rm 128
\yr 1997
\pages 89-157
\endref

\ref \key 8 \by D. Bisch \& V.F.R. Jones \paper \rm A note on free composition of subfactors \inbook  \it
``Geometry and Physics" \publ Marcel Dekker, Lecture Notes in Pure and Applied Mathematics \vol \rm 184 \yr 1997,
339-361
\endref

\ref
\key 9
\by D. Bisch \& V.F.R. Jones
\paper \rm Singly generated planar algebras of small dimension
\jour \it Duke Math. Journal
\vol \rm 101
\yr 2000
\pages 41-75
\endref

\ref
\key 10
\by D. Bisch \& V.F.R. Jones
\paper \rm Singly generated planar algebras of small dimension, Part II
\jour {\it Advances in Math.} (to appear)
\endref

\ref
\key 11
\by D. Bisch \& S. Popa
\paper \rm Examples of subfactors with property T standard invariant
\jour  \it Geom. Funct. Anal.
\vol \rm 9
\yr 1999
\pages 215-225
\endref

\ref
\key 12
\by D. Bisch \& S. Popa
\paper \rm A continuous family of non-isomorphic
irreducible hyperfinite subfactors with the same standard invariant
\jour \it in preparation.
\endref

\ref
\key 13
\by A. Connes
\book \it Noncommutative geometry
\publ Academic Press
\yr 1994
\endref

\ref
\key 14
\by D. Evans \& Y. Kawahigashi
\book \it Quantum symmetries on operator algebras
\publ Oxford University Press
\yr 1998
\endref

\ref
\key 15
\by P. Di Francesco
\paper \rm New integrable lattice models from Fuss-Catalan algebras
\jour \it  Nuclear Phys. B
\vol \rm 532
\yr 1998
\pages 609-634
\endref

\ref
\key 16
\by F. Goodman \& P. de la Harpe \& V.F.R. Jones
\book \it Coxeter graphs and towers of algebras
\publ Springer Verlag, MSRI publications
\yr 1989
\endref

\ref \key 17 \by U. Haagerup \paper \rm Principal graphs of subfactors in the index range $4 < [M:N] < 3 +
\sqrt{2}$ \inbook Subfactors (Kyuzeso, 1993) \publ World Sci. Publishing, River Edge, NJ \yr 1994, 1--38
\endref

\ref
\key 18
\by M. Izumi
\paper \rm Applications of fusion rules to classification of subfactors
\jour \it Publ. RIMS, Kyoto Univ.
\yr 1991
\vol \rm 27
\pages 953-994
\endref

\ref
\key 19
\by V.F.R. Jones
\paper \rm Index for subfactors
\jour \it Invent. Math.
\vol \rm 72
\pages 1-25
\yr 1983
\endref

\ref
\key 20
\by V.F.R. Jones
\paper \rm Hecke algebra representations of braid groups and link polynomials
\jour \it Ann. of Math.
\vol \rm 126
\pages 335-388
\yr
\endref

\ref
\key 21
\by V.F.R. Jones
\paper \rm On a certain value of the Kauffman polynomial
\jour \it Comm. Math. Phys.
\vol \rm 125
\pages 459--467
\yr 1989
\endref

\ref
\key 22
\by V.F.R. Jones
\paper \rm  Planar algebras I
\jour  \it preprint
\endref

\ref \key 23 \by V.F.R. Jones \paper \rm The planar algebra of a bipartite graph \inbook \it Knots in Hellas '98
(Delphi) \publ World Sci. Publishing \yr 2000, 94-117
\endref

\ref
\key 24
\by V.F.R. Jones
\paper \rm The annular structure of subfactors
\jour {\it Enseign. Math.} (to appear)
\endref

\ref
\key 25
\by L. Kauffman
\paper \rm State models and the Jones polynomial
\jour \it Topology
\vol \rm 26
\yr 1987
\pages 395-407
\endref

\ref
\key 26
\by Z. Landau
\paper \rm Fuss-Catalan algebras and chains of intermediate subfactors
\jour \it Pacific J. Math.
\vol  \rm 197
\yr 2001
\pages 325-36
\endref

\ref
\key 27
\by Z. Landau
\paper \rm Exchange relation planar algebras
\jour \it preprint
\yr 2000
\endref

\ref
\key 28
\by J.P. May
\paper \rm Definitions: operads, algebras and modules
\jour \it Contemporary Mathematics
\vol \rm 202
\yr 1997
\pages 1-7
\endref

\ref
\key 29
\by M. J. Martins \& B. Nienhuis
\paper \rm Applications of Temperley-Lieb algebras to Lorentz lattice gases
\jour \it J. Phys. A
\vol \rm 31
\yr 1998
\pages L723--L729
\endref

\ref
\key 30
\by A. Ocneanu
\paper \rm Quantized group string algebras and Galois theory for operator
algebras, in Operator Algebras and Applications 2
\jour \it London Math. Soc. Lect. Notes Series
\vol \rm 136
\yr 1988
\pages 119-172
\endref

\ref
\key 31
\by M. Pimsner \& S. Popa
\paper \rm Entropy and index for subfactors
\jour \it Ann. scient. Ec. Norm. Sup.
\vol \rm 19
\yr 1986
\pages 57-106
\endref

\ref
\key 32
\by S. Popa
\paper \rm Classification of subfactors: reduction to commuting squares
\jour \it Invent. Math.
\vol \rm 101
\yr 1990
\pages 19-43
\endref

\ref
\key 33
\by S. Popa
\paper \rm Classification of amenable subfactors of type II
\jour \it Acta Math.
\vol \rm 172
\yr 1994
\pages 352-445
\endref

\ref
\key 34
\by S. Popa
\paper \rm An axiomatizaton of the lattice of higher relative commutants
\jour \it Invent. Math.
\yr 1995
\vol \rm 120
\pages 427-445
\endref

\ref
\key 35
\by A. Wassermann
\paper \rm Operator algebras and conformal field theory III
\jour \it Invent. Math.
\vol \rm 92
\pages 467-538
\yr 1998
\endref

\ref
\key 36
\by H. Wenzl
\paper \rm Quantum groups and subfactors of type $B$, $C$ and $D$
\jour \it Comm. Math. Phys
\vol \rm 133
\year 1990
\pages 383-432
\endref

\ref
\key 37
\by H. Wenzl
\paper \rm $C\sp *$ tensor categories from quantum groups
\jour \it J. Amer. Math. Soc.
\vol \rm 11
\yr 1998
\pages 261-282
\endref

\ref
\key 38
\by F. Xu
\paper \rm Standard $\lambda$-lattices from quantum groups
\jour \it Invent. Math.
\vol \rm 134
\yr 1998
\pages 455--487
\endref

\enddocument